\documentclass[conference]{IEEEtran}

\usepackage{graphicx}
\usepackage{pstricks}
\usepackage{amssymb}
\usepackage{mathptmx}
\usepackage{amsmath}
\usepackage{amsthm}
\usepackage{MnSymbol}
\usepackage{epstopdf}
\usepackage{balance}
\usepackage{fancyhdr}
\usepackage[colorlinks=true, pdfstartview=FitV, linkcolor=blue, citecolor=blue, urlcolor=blue]{hyperref}

\DeclareMathAlphabet{\mathcal}{OMS}{cmsy}{m}{n}

\renewcommand{\d}{\text{d}}

\newcommand{\pvint}{-\hspace*{-3.5mm}\int}

\newcommand{\phasor}[1]{\mathbf{#1}}

\usepackage{eso-pic}

\begin{document}

\title{Active and Reactive Energy Balance Equations\\ in Active and Reactive Time}

\author{\IEEEauthorblockN{Dimitri Jeltsema}
\IEEEauthorblockA{Delft Institute of Applied Mathematics\\
Delft University of Technology, NL\\
Email: d.jeltsema@tudelft.nl} \and
\IEEEauthorblockN{Gerald Kaiser}
\IEEEauthorblockA{Center for Signals and Waves\\
Portland, OR, USA\\
Email: kaiser@wavelets.com
}}

\maketitle

\begin{abstract}
Electrical networks, and physical systems in general, are known to satisfy a power balance equation which states that the rate of change of the energy in time equals the power at the port of the network minus the power dissipated. However, when complex power is considered, there does not seem to exist a similar statement for the imaginary power, either in the time-domain or in the frequency-domain. Recently, in the context of electromagnetic fields, it has been shown by complexifying the time to \emph{t}\,+\,\emph{js} and interpreting \emph{s} as \emph{reactive time}, that it is possible to set up an imaginary power balance in terms of the rate of change of reactive energy in reactive time. Here these ideas are specialized to linear and time-invariant RLC networks. For non-sinusoidal waveforms it is shown that the rate of change of reactive energy in reactive time contains all the essential properties and features of the commonly accepted definition of reactive power under sinusoidal conditions. We believe that this provides an unambiguous and physically motivated resolution to the longstanding debate on how to generalize reactive power to non-sinusoidal waveforms.  
\end{abstract}

\IEEEpeerreviewmaketitle

\section{Introduction}

The ongoing quest for generalizing the triangular power equality $S^2=P^2+Q^2$ dates back to the work of Steinmetz \cite{Steinmetz1892}.\footnote{Unfortunately, this seminal work is often cited using erroneous page and volume numbers. The correct reference, including a downloadable English translation, is given in \cite{Steinmetz1892}.} In $1892$, he observed that the apparent power $S$ in electric arcs with non-sinusoidal alternating current waveforms can be higher than the active power $P$, even when the reactive power $Q$ appears to be zero. To this day, there is still no consensus about how this difference should be explained under all given load situations. This has resulted in a variety of mutually inconsistent power models.

Traditionally, the problem is either approached in the time domain or in the frequency domain, or sometimes in a combination of both; see for instance \cite{Pasko2012,EmanuelBook} for an overview of the most prominent power models. Apart from the discussion of whether the analysis should be performed in the time or the frequency domains, a new definition of reactive power is usually introduced as follows \cite{Filipski1993}. First, if the waveforms are periodic, the active power is defined as the periodic mean of the instantaneous power, i.e.,
\begin{equation}\label{eq:Pactive}
P := \frac{1}{T} \int\limits_0^T u(t)i(t)\d t = \langle u,i \rangle
\end{equation} 
and the apparent power as $S=\|u\|\,\|i \|$, with $\| \cdot \| = \sqrt{\langle \, \cdot \, , \, \cdot \,\rangle}$ the rms norm. Then, after introducing a new definition of the reactive power, say $Q_\text{new}$, it is assumed that the geometric sum of all power components is equal to the apparent power. When the geometric sum is smaller than the apparent power, a new quantity, say REST, is introduced to fulfil 
\begin{equation}\label{eq:S + rest}
S^2 = P^2 + Q_\text{new}^2 + (\text{REST})^2.
\end{equation} 
Only in Fryze's (time-domain) power model is the `REST' quantity absent. In most other theories `REST' is either considered as an essential measure of a certain load parameter or is decomposed further into the geometric sum of new power components, $(\text{REST})^2 = (\text{REST}_1)^2 + (\text{REST}_2)^2 + \cdots$. 

\subsection{Motivation for a Different Approach}

Power models along the line of (\ref{eq:S + rest}) have led to useful models for compensator design, and to a certain extent they have improved our understanding of what is going on. Perhaps the most detailed example is the \emph{Currents' Physical Components} (CPC) model by Czarnecki \cite{Czarnecki1991,Czarnecki2008} (with over $5$ `REST' components)  or the $7/11$ component decomposition  of \cite{Hanoch2006}. Our main objection, however, is that such models have no solid physical foundation in terms of conserved quantities, i.e., no energy and power balance equations. Additionally, some of the power components of the aforementioned decompositions are not unique in the sense that they depend on the chosen reference. For example, using the voltage as a reference generally leads to a different value of reactive power than using the current as a reference \cite{Czarnecki1991}.

In \cite{JeltsemaPrzeglad2016}, the first author of the present paper went back to the classical sinusoidal phasor-based complex power model and generalized the phasor concept to periodic non-sinusoidal waveforms using the Hilbert transform. This approach provides a consistent generalization of the classical power trinity to its time-varying counterpart $\{P(t),Q(t),S(t)\}$, referred to as the real power, the imaginary power, and the time-varying apparent power, respectively. The real power $P(t)$ can be associated to a power balance which states that the rate of change of a potential function that is related to the total stored magnetic and electric {energies} in time equals the power at the network port minus the dissipated power. The time-average of this balance provides the active power as defined in (\ref{eq:Pactive}). On the other hand, for linear and time-invariant (LTI) systems $Q(t)$ is \emph{related} to the difference between the total stored magnetic and electric energy, and its mean value provides Budeanu's reactive power \cite{Budeanu1927}, but it cannot be associated with the rate of change of some potential function in time, to be interpreted as reactive energy. From a physical perspective, this is rather dissatisfying.

\subsection{The Generalized Complex Poynting Theorem}

Recently, in the context of general electromagnetic fields in free space, the second author of the present paper has derived a \emph{generalized complex Poynting theorem} (GCPT) in which the active and reactive energies\footnote{The term \emph{active} is used here since active energy can perform work, while \emph{reactive} energy cannot.}  play completely symmetrical roles \cite{Kaiser}. This is accomplished by replacing all waveforms by their \emph{analytic signals} (positive-frequency parts) and then complexifying the time, $t\to \tau=t+js$, with $s>0$ and $j^2=-1$, so that $e^{j\omega\tau}$ gives absolute convergence of all inverse Fourier integrals to analytic functions of $\tau$. It was then shown that (a) $s$ is a \emph{time-resolution scale} and (b) all analytic signals are moving averages of the original signals over a `fuzzy' moving time interval $\approx(t-s, t+s)$. The real part of the GCPT then gives the power balance for the rate of change of  the \emph{active energy} in $t$,  while its imaginary part gives the power balance for the rate of change of \emph{reactive energy} in $s$ (not $t$). All energies and powers are moving averages, as explained above. For this reason, $t$ is now reinterpreted as the \emph{active time (interval)}  and $s$ as the \emph{reactive time}, the latter to be measured in \emph{seconds reactive} [sr]. The GCPT thus \emph{completes} the time-harmonic complex Poynting theorem by giving an \emph{equal but dual} physical status to active and reactive energies and powers.\footnote{For $s=0$, the GCPT reduces to the \emph{real} Poynting theorem, thus bridging the gap between the real and complex Poynting theorems.} 
Furthermore, it provides a rational basis for the Volt-Amp\`ere reactive [VAr], the somewhat strange unit used to measure reactive power. Namely, the reactive current in the balance law for the reactive power is the flow of charge with respect to $s$, measured in Coulombs per seconds reactive [C/sr] or Amp\`eres reactive [Ar], and the reactive power is the rate of change of the reactive energy with respect to $s$, measured in Joules per second reactive\,\footnote{Reactive energy, like active energy, is measured in Joules. Only the units for \it rates of change \rm are different for active and reactive quantities, depending on whether the rates are with respect to $t$ or $s$.} or Volt-Amp\`ere 
reactive:
\begin{align*}\label{VAr} 
\boxed{\rm [J/sr]=[V\times C/sr]= [VAr].}
\end{align*}
Without the notion of reactive time, the [VAr] seems \emph{ad hoc}.

\subsection{Contribution and Outline}

The present paper is a circuit-theoretic counterpart of \cite{Kaiser} and a refinement of \cite{JeltsemaPrzeglad2016}. The contribution here is two-fold: First, in Sections \ref{sec:balances}--\ref{sec:classical}, the necessary background material is provided on balance equations and the classical sinusoidal power model. It is pointed out that the imaginary part of the classical sinusoidal complex power model does not allow for a power balance equation. Section \ref{sec:REints} extends the classical phasor representation to general (not necessarily sinusoidal or even periodic) waveforms using the \emph{analytic-signal transform}. This makes possible the symmetry between the real and imaginary parts of the complex power and finally provides a balance equation for each part. Secondly, Section \ref{sec:Xinstant} shows that there cannot exist an `instantaneous' reactive energy and refines the ideas of \cite{JeltsemaPrzeglad2016} concerning Budeanu's reactive power. The results are illustrated using a linear and time-invariant (LTI) RLC load subject to a sinusoidal voltage source that contains a flicker component. In Section \ref{sec:s} we offer some final remarks.

\section{Power and Energy Balance Equations}\label{sec:balances}

Consider a single-phase generator transmitting energy to a load. We know from Tellegen's theorem that at each instant the rate of energy (i.e., the power) entering the load through its port gets distributed among the elements of the load network, so that none is lost. For a single-phase load network consisting of branches $b$ (edges), the instantaneous power at the port is 
\begin{equation}\label{eq:PB}
p(t):=u(t)i(t) = \sum_b u_b(t)i_b(t).
\end{equation}
This equation holds regardless of the nature of the elements or the excitation and is consistent with the principle of conservation of energy. If the branches consists of LTI resistors, inductors, and capacitors, then (\ref{eq:PB}) can be explicitly written as
\begin{equation*}
\begin{aligned}
p(t) = \sum_{\text{res}} R_b\, i_b^2(t) + \sum_{\text{ind}} L_b \frac{\d i_b}{\d t}(t)i_b(t) + \sum_{\text{cap}}  u_b(t) C_b\frac{\d u_b}{\d t}(t),
\end{aligned}
\end{equation*}
or, equivalently, as
\begin{equation}\label{eq:PBE}
\boxed{ \frac{\d w}{\d t}(t) =  p(t) - p_d(t),}
\end{equation}
where $w(t) := w_m(t) + w_e(t)$ is the total instant magnetic and electric energy, with
\begin{equation}\label{eq:wme}
w_m(t) = \frac{1}{2}\sum_{\text{ind}} L_b i_b^2(t), \quad w_e(t) = \frac{1}{2}\sum_{\text{cap}} C_b u_b^2(t)
\end{equation}
and
\begin{equation}\label{eq:pd}
p_d(t) = \sum_{\text{res}} R_b i_b^2(t)
\end{equation}
the total instantaneous dissipated power.

Integrating (\ref{eq:PBE}) over some time-interval, say $[t_0,t_1]$, gives
\begin{equation}\label{eq:EB}
\boxed{w(t_1) - w(t_0) = \int\limits_{t_0}^{t_1} p(t) \d t - \int\limits_{t_0}^{t_1} p_d(t) \d t,}
\end{equation}
which states that the total stored energy equals the energy at the port minus the total dissipated energy. For that reason, we refer to (\ref{eq:EB}) as an \emph{energy balance} equation and to (\ref{eq:PBE}) as a \emph{power balance} equation. The instantaneous (or \emph{local}) balance (\ref{eq:PBE}) thus implies the `global' balance (\ref{eq:EB}). 

\section{The Classical Sinusoidal Power Model}\label{sec:classical}

Let the voltage at the load terminals be given by
\begin{equation}\label{eq:u_sinus}
u(t)=U\sqrt{2}\cos(\omega t + \alpha),
\end{equation}
where $\omega = 2 \pi/T$ and $U$ is the rms voltage.
Under the assumption that $u(t)$ does not depend on the transmitted current (infinitely strong power system), the associated current reads
\begin{equation}\label{eq:i_sinus}
i(t)=I\sqrt{2}\cos(\omega t + \beta).
\end{equation}
The instantaneous power can be shown to be
\begin{equation}\label{eq:p_sinus}
p(t) = P\big[1+\cos(2\omega t +2\alpha)\big] + Q\sin(2\omega t+2\alpha),
\end{equation}
where $P$ and $Q$ are the active and reactive powers defined by
\begin{equation}\label{eq:PandQsinus}
\begin{aligned} 
P:=UI\cos(\varphi),\quad
Q:=UI\sin(\varphi),
\end{aligned}
\end{equation}
with $\varphi:=\alpha - \beta$ the phase-shift between $u(t)$ and $i(t)$.

Since $w(T)=w(0)$, Equation (\ref{eq:EB}) implies that
\begin{equation}\label{eq:Pdaverage}
P= \frac{1}{T}\int\limits_0^T p_d(t) \d t,
\end{equation}
which states that the \emph{active} power is the time average of the power dissipated in a period.

Unfortunately, the energy balance equation (\ref{eq:EB}) reveals little about the physical nature of  the second term in (\ref{eq:p_sinus}). There is no obvious connection between (\ref{eq:PBE}) and the reactive power $Q$ similar to (\ref{eq:Pdaverage}).

\subsection{Phasors and Complex Power}

Alternatively, a standard method in electrical engineering is to represent (\ref{eq:u_sinus}) and (\ref{eq:i_sinus}) by their time-harmonic \emph{phasors}  \cite{Desoer}
\begin{equation}\label{eq:ui_rotating}
\phasor{u}(t) = \phasor{U}e^{j\omega t}, \quad \phasor{i}(t) = \phasor{I}e^{j\omega t},
\end{equation}
with $\phasor{U}=U\sqrt{2}e^{j\alpha}$ and $\phasor{I}=I\sqrt{2}e^{j\beta}$. This enables one to define the \emph{complex power}
\begin{equation}\label{eq:Scomplex}
\phasor{S}:=\frac{1}{2}\phasor{u}(t)\phasor{i}^*(t) = UIe^{j\varphi} = P + jQ,
\end{equation}
with its magnitude equal to the apparent power, i.e., $\left|\,\phasor{S}\,\right|=S$. 

Furthermore, for arbitrary LTI RLC loads, the phasor counterpart of (\ref{eq:PB}) is given by 
\begin{equation}\label{eq:TellegenComplexRLC}
\begin{aligned}
\frac{1}{2}\sum_b &\phasor{u}_b(t)\phasor{i}_b^*(t) = \frac{1}{2}\sum_{\text{res}} R_b\,\phasor{i}_b(t)\phasor{i}_b^*(t)\\
& \quad  + \frac{1}{2}\sum_{\text{ind}} L_b\frac{\d\phasor{i}_b}{\d t}(t)\phasor{i}_b^*(t) + \frac{1}{2}\sum_{\text{cap}} \phasor{u}_b(t)C_b\frac{\d\phasor{u}_b^*}{\d t}(t),
\end{aligned}
\end{equation}
which, upon substitution of $\phasor{u}_b(t)=\phasor{U}_be^{j\omega t}$ and $\phasor{i}_b(t)=\phasor{I}_be^{j\omega t}$, yields the well-know expression \cite{Desoer}
\begin{equation}\label{eq:S_RLC}
\boxed{ 
\phasor{S} := \frac{1}{2}\sum_\text{res} R_b \left|\,\phasor{I}_b\,\right|^2 +j2\omega \left[ \frac{1}{4}\sum_\text{ind}L_b \left|\,\phasor{I}_b\,\right|^2 - \frac{1}{4}\sum_\text{cap}C_b \left|\,\phasor{U}_b\,\right|^2 \right].
}
\end{equation}
Comparing the latter with (\ref{eq:Scomplex}) reveals that $P=P_d$ with
\begin{equation}\label{eq:Pd_sinus}
P_d:=\frac{1}{2}\sum_\text{res} R_b \left|\,\phasor{I}_b\,\right|^2,
\end{equation}
whereas, setting
\begin{equation}\label{eq:Wm_We_sinus}
W_m := \frac{1}{4}\sum_\text{ind}L_b \left|\,\phasor{I}_b\,\right|^2, \quad
 W_e:=\frac{1}{4}\sum_\text{cap}C_b \left|\,\phasor{U}_b\,\right|^2,
\end{equation}
yields  
\begin{equation}\label{eq:QWmWe}
Q= 2\omega \left[ W_m - W_e \right],
\end{equation}
where $W_m$ and $W_e$ are the mean values of $w_m(t)$ and $w_e(t)$ as given in (\ref{eq:wme}).

\subsection{Existence of Reactive Energy in Active Time}\label{sec:TVphasor}

Since $\d W / \d t =0$ as $W:=W_m+W_e$ is constant, the real part of (\ref{eq:S_RLC}) gives $P=P_d$ (which also directly follows from (\ref{eq:Pdaverage})). However, there is no such balance equation for the imaginary part of (\ref{eq:S_RLC}) in terms of a (rate of change of) reactive energy in active (real) time. The obstacle here is that the imaginary part is not integrable with respect to the voltage and current phasors, and therefore does not allow for a potential function that can be associated to reactive energy. This fact  was proved in \cite{Carozzi2005} in the context of electromagnetic fields.

\section{Reactive Energy in the Time-Scale Domain}\label{sec:REints}

To resolve the asymmetry between the real and imaginary parts of the complex power, and at the same time providing a fully consistent generalization of the classical power model for \emph{arbitrary} square-integrable waveforms, we proceed along the lines of \cite{Kaiser}.  

\subsection{Phasors and Hilbert Transform}

\begin{figure}[t]
\begin{center}
\includegraphics[width=80mm]{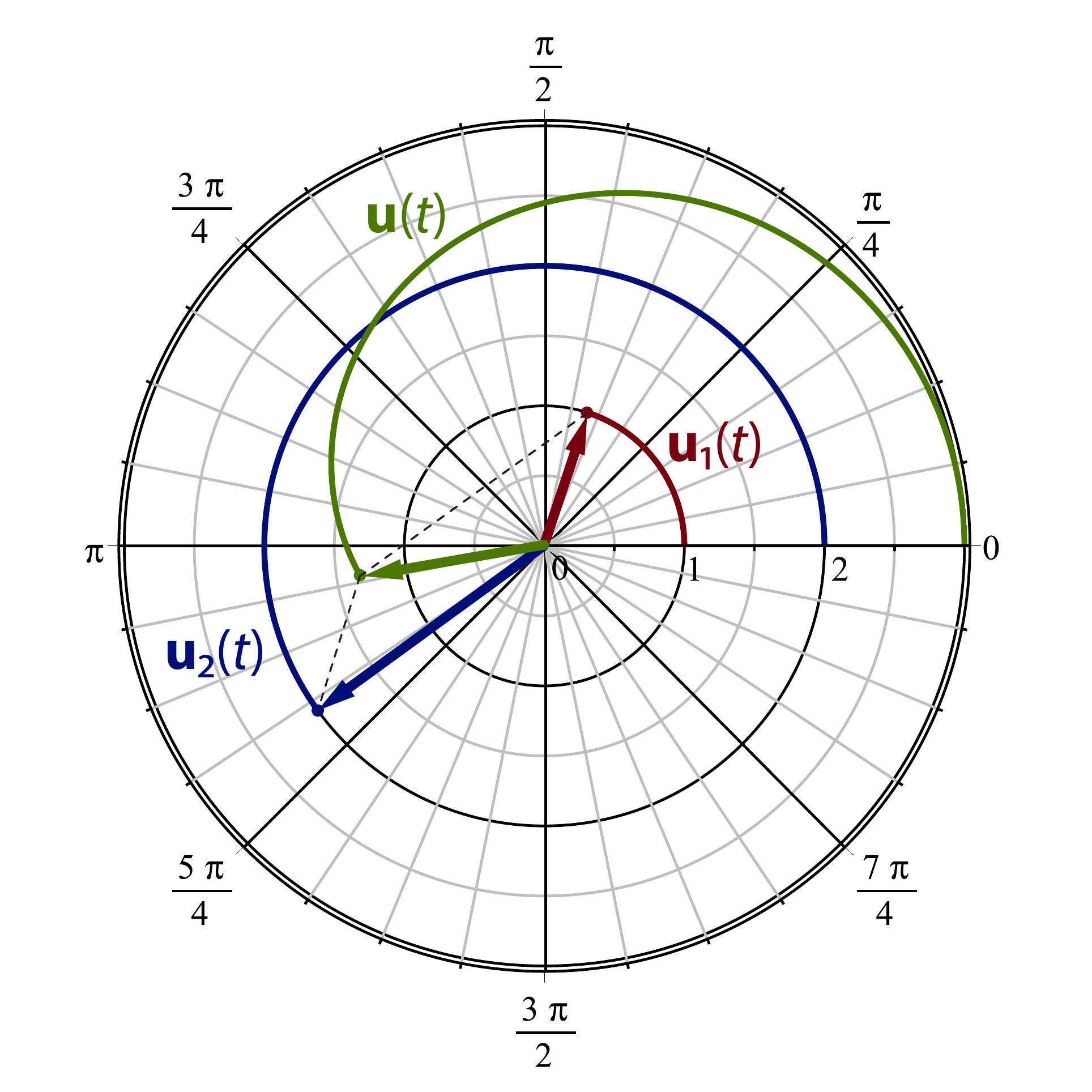}
\caption{The analytic signals (phasors) associated to $u_1(t) = \cos(t)$, $u_2(t) =2\cos(3t)$, and $u(t)= u_1(t) + u_2(t)$, for $0 \leq t \leq 2\pi/5$. Observe that $\phasor{u}_1(t)=e^{\,jt}$ rotates at angular velocity $1$ rad/s and  $\phasor{u}_2(t)=2e^{\, j3t}$ at $3$ rad/s, with magnitudes $1$ and $2$, respectively. At each time instant, the analytic voltage $\phasor{u}(t)$ obtained from (\ref{eq:u_analytic}) represents the (vector) sum of $\phasor{u}_1(t)$ and $\phasor{u}_2(t)$.}
\label{fig:phasors}
\end{center}
\end{figure}

The underlying mathematical principle behind the transition from the sinusoidal time functions (\ref{eq:u_sinus})--(\ref{eq:i_sinus}) to their rotating phasor representation (\ref{eq:ui_rotating}) is the \emph{analytic signal}, widely used in telecommunication and signal processing \cite{VakmanBook}. For an arbitrary real voltage-current pair $f(t)=\{u(t),i(t)\}$, the analytic signal representation is defined as the \emph{positive-frequency} part of the inverse Fourier transform:
\begin{equation}\label{eq:u_analytic}
\phasor{f}(t) = \frac{1}{\pi} \int\limits_{0}^\infty \hat{F}(\omega)e^{j\omega t} \d \omega = f(t) + j f_h(t), 
\end{equation}  
where $\hat{F}(\omega)$ denotes the frequency spectrum and $f_h(t)$ denotes the \it Hilbert transform \rm of $f(t)$, given by
\begin{equation}
f_h(t) = \frac{1}{\pi} \pvint\limits_{-\infty}^\infty \frac{f(t')}{t-t'} \d t',
\end{equation}
where the Cauchy principal value of the integral is implied. The geometric interpretation of  (\ref{eq:u_analytic}) is illustrated in Fig.~\ref{fig:phasors} for a voltage composed of two sinusoidal components.

\subsection{Analytic-Phasor Transform}

Since the integral in (\ref{eq:u_analytic}) is restricted to $\omega > 0$, it extends analytically to the upper half of the complex time plane by 
\begin{equation}\label{eq:u_analytic_time}
\phasor{f}(t+js) = \frac{1}{\pi} \int\limits_{0}^\infty \hat{F}(\omega)e^{j\omega(t+js)} \d \omega, \ s>0. 
\end{equation}
It is shown in \cite{Kaiser} that $s$ represents a \emph{time-resolution scale},  so we call $(t,s)$ the \emph{time-scale} domain.\footnote{Since the transformation (\ref{eq:u_analytic_time}) provides windowed time averages over the interval $\Delta\, t  \sim \pm s$, the time-scaled waveforms are \emph{non-local} in character. This is the basis for interpreting $s$ as a \emph{time-resolution scale}. Note that in the limit $s \to 0$, we have $\phasor{f}(t+js) \to \phasor{f}(t)$; see \cite{Kaiser} for more details.}
The integrand of (\ref{eq:u_analytic_time}) gains a low-pass filter with transfer function $e^{-\omega s}$, which suppresses frequencies $\omega \gg 1/s$. Furthermore, if $f(t)$ is square-integrable (`finite-energy'), then (\ref{eq:u_analytic_time}) defines an \emph{analytic} voltage-current pair $\phasor{f}(t+js)$ in the upper-half complex time plane.

A direct time-domain representation of (\ref{eq:u_analytic_time}) is given by the \emph{analytic-phasor transform} \cite{Kaiser}
\begin{equation}\label{eq:ast}
\boxed{\phasor{f}(t+js) = \frac{j}{\pi} \int\limits_{-\infty}^\infty \frac{f(t')}{t+js - t'}\d t'.} 
\end{equation}
Roughly, $\phasor{f}(t+js)$ depends on $f(t')$ mainly in the interval $t \pm s$. Along with the property
\begin{align*}
\lim_{s\to 0}\, \text{Re}\{\phasor{f}(t+js)\} =f(t),
\end{align*}
this justifies calling $s$ a time-resolution scale. 

\subsection{Time-Scaled Complex Power Balance}

Let $\tau:= t+js$, then (\ref{eq:TellegenComplexRLC}) in the time-scale $(t,s)$ extends to
\begin{equation*}
\begin{aligned}
\frac{1}{2} &\phasor{u}(\tau)(\phasor{i}(\tau))^* = \frac{1}{2}\sum_{\text{res}} R_b\phasor{i}_b(\tau)(\phasor{i}_b(\tau))^*\\
& \quad  + \frac{1}{2}\sum_{\text{ind}} L_b\frac{\partial\phasor{i}_b(\tau)}{\partial\tau}(\phasor{i}_b(\tau))^* + \frac{1}{2}\sum_{\text{cap}} \phasor{u}_b(\tau)C_b\frac{\partial(\phasor{u}_b(\tau))^*}{\partial\tau^*}.
\end{aligned}
\end{equation*}
Since $\frac{\partial}{\partial\tau}(\phasor{i}_b(\tau))^*=0$ and $\frac{\partial}{\partial\tau^*}\phasor{u}_b(\tau)=0$ by analyticity,\footnote{Since $t=\frac{1}{2}(\tau + \tau^*)$ and $s=j\frac{1}{2}(\tau^* - \tau)$, it follows that
\vspace*{-0.5em}
\begin{align*}
\tfrac{\partial}{\partial\tau}= \tfrac{1}{2}\left[\tfrac{\partial}{\partial t} - j \tfrac{\partial } {\partial s} \right], \quad
\tfrac{\partial}{\partial\tau^*} = \tfrac{1}{2}\left[\tfrac{\partial}{\partial t} + j \tfrac{\partial } {\partial s} \right].
\end{align*}
} this is equivalent to
\begin{equation}\label{eq:TellegenComplexRLC_ts}
\begin{aligned}
\frac{1}{2} &\phasor{u}(\tau)(\phasor{i}(\tau))^* = \frac{1}{2}\sum_{\text{res}} R_b|\phasor{i}_b(\tau)|^2\\
& \quad  + \frac{1}{2}\frac{\partial}{\partial\tau}\sum_{\text{ind}} L_b|\phasor{i}_b(\tau)|^2 + \frac{1}{2}\frac{\partial}{\partial\tau^*}\sum_{\text{cap}} C_b|\phasor{u}_b(\tau)|^2.
\end{aligned}
\end{equation}
The real and imaginary parts are now fully symmetric. To see this, define the \emph{scaled} magnetic and electric energies 
\begin{align*}
\mathcal{W}_m(t,s) := \frac{1}{4}\sum_{\text{ind}}L_b\left|\phasor{i}_b(\tau)\right|^2, \quad  
\mathcal{W}_e(t,s) := \frac{1}{4}\sum_{\text{cap}}C_b\left|\phasor{u}_b(\tau)\right|^2,
\end{align*}
define the \emph{scaled} dissipated power
\begin{equation*}
\mathcal{P}_d(t,s) := \frac{1}{2}\sum_{\text{res}} R_b |\phasor{i}_b(\tau)|^2,
\end{equation*}
and the \emph{scaled} real and imaginary powers
\begin{align*}
\mathcal{P}(t,s) =\frac{1}{2}\text{Re}\big\{\phasor{u}(\tau)(\phasor{i}(\tau))^*\big\}, \quad
\mathcal{Q}(t,s) =\frac{1}{2}\text{Im}\big\{\phasor{u}(\tau)(\phasor{i}(\tau))^*\big\}.
\end{align*}
Then (\ref{eq:TellegenComplexRLC_ts}) in the time-scale domain becomes
\begin{equation}\label{eq:TellegenComplexRLC_ts_final}
\boxed{\frac{\partial \mathcal{W}}{\partial t}(t,s) - j \frac{\partial \mathcal{X}}{\partial s}(t,s) = \mathcal{P}(t,s) - \mathcal{P}_d(t,s)  + j \mathcal{Q}(t,s),}
\end{equation}
where $\mathcal{W}  := \mathcal{W}_m + \mathcal{W}_e$ and $\mathcal{X}  := \mathcal{W}_m - \mathcal{W}_e$ represent the \emph{scaled} active and reactive energy, respectively. 

Note that (\ref{eq:TellegenComplexRLC_ts_final}) naturally splits into an active and a reactive power balance law, i.e.,
\begin{equation}\label{eq:PBactive_st}
\boxed{\frac{\partial \mathcal{W}}{\partial t}(t,s) = \mathcal{P}(t,s) - \mathcal{P}_d(t,s),}
\end{equation}
and
\begin{equation}\label{eq:PBreactive_st}
 \boxed{-\frac{\partial \mathcal{X}}{\partial s}(t,s) = \mathcal{Q}(t,s),}
\end{equation}
respectively. 

Clearly, (\ref{eq:PBactive_st}) constitutes a power balance equation in the active time $t$ for the rate of change of the average active energy at scale $s$. This is the scaled time-averaged counterpart of (\ref{eq:PBE}). The imaginary part (\ref{eq:PBreactive_st}), however, does not constitute a power balance equation in the active time $t$, but in the scale $s$. It states that scaled reactive power equals the rate of change of scaled reactive energy with respect to  \emph{scale refinements}.\footnote{Due to the sign of $-\frac{\partial}{\partial s}\mathcal{X}(t,s)$ and the fact that $\mathcal{X}(t,s) \to 0$ as $s \to \infty$, the orientation of $s$ is from \emph{course} to \emph{fine} scales and $\mathcal{X}(t,s)$ represents the \emph{cumulative} reactive energy at all scales $s' \geq s$; see \cite{Kaiser}.}  This may seem rather unconventional, but we simply have no choice, as demonstrated in Section \ref{sec:TVphasor}. It is simply the price of having a symmetric description of active and reactive energy.

\subsection{The Sinusoidal Case}

At this point it is instructive to examine how for sinusoidal waveforms (\ref{eq:TellegenComplexRLC_ts_final}) relates to (\ref{eq:S_RLC}). Since $\phasor{u}(\tau)=\phasor{U}e^{j\omega t}e^{-\omega s}$ and $\phasor{i}(\tau)=\phasor{I}e^{j\omega t}e^{-\omega s}$, we have 
\begin{equation*}
\left|\,\phasor{u}(\tau)\,\right|^2 = e^{-2\omega s}\left|\,\phasor{U}\,\right|^2, \ \left|\,\phasor{i}(\tau)\,\right|^2 = e^{-2\omega s}\left|\,\phasor{I}\,\right|^2,
\end{equation*}
which clearly are time-independent. Consequently, all scaled quantities in   (\ref{eq:TellegenComplexRLC_ts_final}) are time-independent, i.e.,
\begin{equation*}
\frac{\partial \mathcal{W}}{\partial t}(s) - j \frac{\partial \mathcal{X}}{\partial s}(s) = \frac{1}{2}e^{-2\omega s} \phasor{U}\,\phasor{I}^* - \frac{1}{2}e^{-2\omega s}\sum_\text{res}R_b \left|\,\phasor{I}_b\,\right|^2.
\end{equation*}
Since $\frac{\partial}{\partial t}\mathcal{W}(t,s)=0$ and $- j \frac{\partial \mathcal{X}}{\partial s}(s) = j2\omega \mathcal{X}(s)$, this reduces to 
\begin{equation}\label{eq:PBs}
\begin{aligned}
\boxed{\phasor{S} e^{-2\omega s} = P_d e^{-2\omega s} + j2\omega\left[W_m - W_e\right]e^{-2\omega s},}
\end{aligned}
\end{equation}
with $P_d$ as defined in (\ref{eq:Pd_sinus}) and $W_m$ and $W_e$ as defined in (\ref{eq:Wm_We_sinus}). Cancelling the $e^{-2\omega s}$ factor gives (\ref{eq:S_RLC}). 

\section{There is no `Instantaneous' Reactive Energy!}\label{sec:Xinstant}

The quantities $\mathcal{W}(t,s)$, $\mathcal{P}(t,s)$, and $\mathcal{P}_d(t,s)$ in (\ref{eq:TellegenComplexRLC_ts_final}) are the scaled versions of the instantaneous  quantities $w(t)$, $p(t)$, and $p_d(t)$ in (\ref{eq:PBE}). Furthermore, $\mathcal{X}(t,s)$ is the scaled version of 
\begin{equation*}
x(t):=w_m(t)-w_e(t),
\end{equation*}
which it is tempting to call the `instantaneous reactive energy'. But here we encounter a fundamental difference between the active and reactive energies, as explained in the more general context of electromagnetic fields in \cite{Kaiser}. 

This difference already becomes  apparent at $s=0$.  
Indeed, for $s \to 0$, the transformation (\ref{eq:u_analytic_time}) reduces to (\ref{eq:u_analytic}), which states that each analytic voltage-current pair $\{\phasor{u}(t),\phasor{i}(t)\}$ is the sum of the \emph{local} (real) part $\{u(t),i(t)\}$ and the \emph{nonlocal} (imaginary) part $\{u_h(t),i_h(t)\}$. In that same limit, $\mathcal{W}(t,s) \to W(t)$, $\mathcal{P}(t,s) \to P(t)$, and $\mathcal{P}_d(t,s) \to P_d(t)$, where 
\begin{equation}\label{eq:P(t)}
P(t) = \frac{1}{2}[u(t)i(t) + u_h(t)i_h(t)] 
\end{equation}
coincides with the \emph{real power} as proposed in \cite{JeltsemaPrzeglad2016}. When these quantities are inserted into (\ref{eq:PBactive_st}), we find that each side of the balance splits into a local part (involving no Hilbert transforms) and a nonlocal part (with Hilbert transforms in each term). Consequently, the local and nonlocal parts \emph{each} satisfy (\ref{eq:PBactive_st}), and the balance equation for the local part simply coincides with the original balance law (\ref{eq:PBE}) for the instantaneous energy. Thus, in the limit $s \to 0$, \emph{active energy can be localized}. 

On the other hand, evaluation of the reactive-energy balance law (\ref{eq:PBreactive_st}) at $s=0$ gives 
\begin{equation}\label{eq:dsW}
\left. -\frac{\partial \mathcal{X}}{\partial s}(t,s)\right|_{s=0} = \mathcal{Q}(t,0),
\end{equation}
with $\mathcal{Q}(t,0) \equiv Q(t)$, and
\begin{equation}\label{eq:Q(t)}
Q(t) = \frac{1}{2}[u_h(t)i(t) - u(t)i_h(t)]
\end{equation}
coincides with the \emph{imaginary power} as proposed in \cite{JeltsemaPrzeglad2016}. Although the reactive energy $\mathcal{X}(t,s)$ can be split into the sum of local and nonlocal parts, these cannot be separated (as in the active case) because they do not individually satisfy (\ref{eq:dsW}). Moreover, for $s>0$, the local and nonlocal parts get mixed together as the scale changes according to (\ref{eq:PBreactive_st}). This shows that \emph{reactive energy is inherently nonlocal}, which suggests that theories of `instantaneous' reactive energy and power must be fundamentally incorrect. 

Thus, although $\mathcal X(t,s), \mathcal Q(t,s)$ and their $s=0$ limits are generally time-dependent, their time variables represent  `fuzzy' time, \emph{not} sharp instants. This may explain why there has been so much controversy and confusion around the origin and generalization of reactive energy and power. 

\subsection{Budeanu's Reactive Power}

In the limit $s \to 0$, the scaled reactive power $\mathcal{Q}(t,s)$ thus reduces to the imaginary power $Q(t)$. Now, since we know from \cite{JeltsemaPrzeglad2016} that Budeanu's reactive power $Q_B$ corresponds to the mean value of (\ref{eq:Q(t)}), it is easily deduced that for arbitrary waveforms $Q_B = \overline{Q(t)}$, with
\begin{equation}\label{eq:Budeanu_st}
\overline{Q(t)} = -\frac{\partial}{\partial s} \left[\lim_{T\to\infty}\frac{1}{T}\int\limits_0^T \big(\mathcal{W}_m(t,s) - \mathcal{W}_e(t,s)\big)\, \d t \right]_{s=0},  
\end{equation}
which states that Budeanu's reactive power equals the \emph{mean} rate of change of the scaled reactive energy at $s=0$. Budeanu's reactive power has the merit to indicate weather the load is dominantly inductive ($Q_B >0$) or  capacitive ($Q_B < 0$). It is indeed not suitable for compensation purposes as properly forcing $Q_B \equiv 0$ only compensates the mean value of the reactive power. Consequently, the load does not exhibit a dominantly inductive or capacitive behaviour anymore, but there still might be a fluctuating part of the reactive power present.  However, in contrast to the assertions made in \cite{Czarnecki1987} against Budeanu's power model, accepting the time-scale domain (\ref{eq:Budeanu_st}) can be considered as proof that Budeanu's reactive power \emph{does} possess a physical interpretation. A resolution of Budeanu's distortion power is proposed in \cite{JeltsemaPrzeglad2016}.

\subsection{Non-sinusoidal Example}

\begin{figure*}[t]
\begin{center}
\includegraphics[width=170mm]{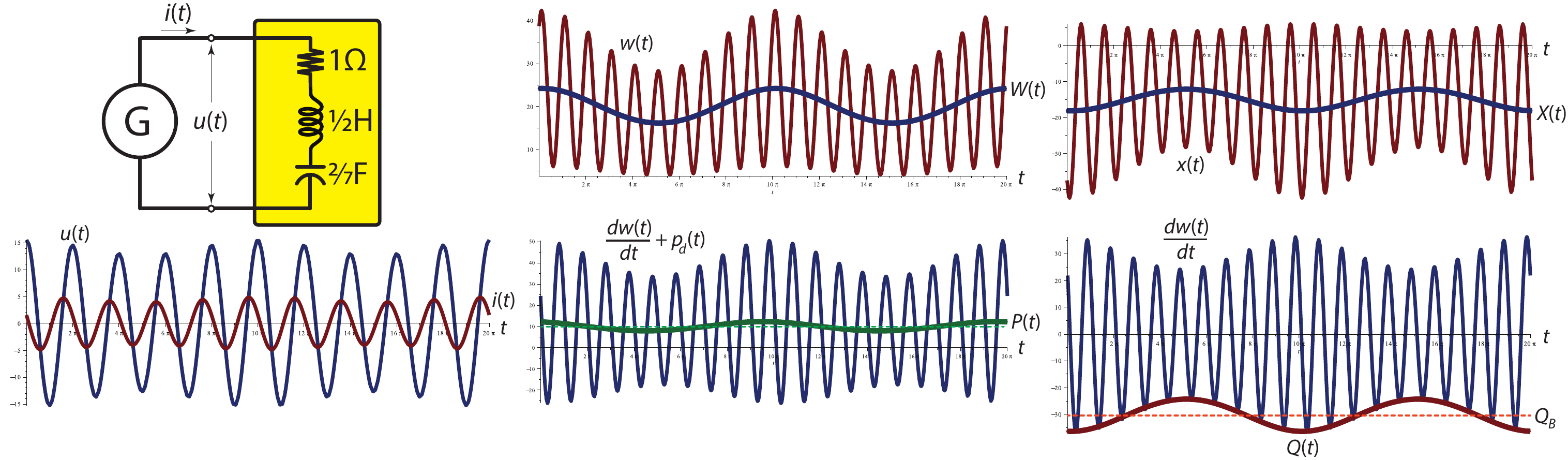}
\caption{LTI RLC load subject to $u(t)=10\sqrt{2}\big(1+0.1\cos(0.2t)\big)\cos(t)$. It is clear from the waveforms that $W(t)$, $X(t)$, and $P(t)$ are moving time averages of $w(t)$, $x(t)$, and $p(t)$, respectively, whereas $Q(t)$ coincides with the time-averaged lower envelope of $\frac{\d}{\d t}w(t)$ since the load is dominantly capacitive. This is also observed from Budeanu's reactive power which follows from taking the mean value of $Q(t)$.}
\label{fig:RLC_flick}
\end{center}
\end{figure*}

Consider the RLC load network of Fig.~\ref{fig:RLC_flick} supplied with a sinusoidal voltage that contains a $10 \%$ flicker component. The resulting energies and powers all become time-varying as shown in Fig.~\ref{fig:RLC_flick} for $s=0$, except for the mean value of (\ref{eq:P(t)}), which represents the active power (\ref{eq:Pactive}), here obtained as
\begin{equation*}
\overline{P(t)}= \lim_{T\to\infty}\frac{1}{T} \int\limits_0^T \left[\frac{\partial \mathcal{W}}{\partial t}(t,s) + \mathcal{P}_d(t,s)\right]_{s=0} \d t =  10.05 \ \text{W},
\end{equation*}
and Budeanu's reactive power $Q_B = \overline{Q(t)}= -30.15$ VAr. 

\section{Concluding Remarks}\label{sec:s}

The analysis in this paper shows that there cannot exist a balance equation for reactive energy in active (real) time. To obtain symmetry between the real and imaginary parts of the complex power (\ref{eq:TellegenComplexRLC}), allowing to  set up  power balance equations for both parts, we must complexify time as $t+js$. This results in a power balance associated to the rate of active energy with respect to $t$ and a power balance associated to the rate of reactive energy with respect to $s$. 

\subsection{Why not Complex Energy?}

It becomes apparent from (\ref{eq:TellegenComplexRLC_ts_final}) that there does not exists an overall scaled `complex potential function' that represents `complex energy' analogous to $\mathcal{P} + j \mathcal{Q}$. There is no \emph{qualitative} difference between active and reactive energy, as both quantities are measured in Joules [J]. The difference is solely \emph{quantitative}:
$\mathcal{W}_e$ contributes positively to $\mathcal{W}$ and negatively to $\mathcal{X}$, whereas $\mathcal{W}_m$ contributes positively to both. 
 Only for rates of change of energy does a qualitative difference emerge between active and reactive variables, as it becomes necessary to specify if the rate of change is with respect to $t$ or $s$. Since $s$ is imaginary time, this also explains why the reactive rates of change enter as imaginary variables in (\ref{eq:TellegenComplexRLC_ts_final}).\footnote{However, since $(js)^2=-s^2$ is real, it follows that $[\rm sr]^2=[s]^2$ and $\rm [VAr]^2=[\rm VA]^2=[\rm W]^2$. This explains why the square of the apparent power $S$ can be computed as $S^2=P^2+Q^2$ and measured in $\rm[W]^2$. From this point of view, it is unnecessary to reserve the units Volt\! $\times$\! Amp\`ere for $S$, as done in \cite{Steinmetz1892}.}

\subsection{Nonlinear and Multi-Phase Load Networks}

The proposed power framework is set up for single-phase LTI load networks. Although the results can be carried over to multi-phase LTI loads replacing (\ref{eq:PB}) with
\begin{equation}\label{eq:PB_poly}
p(t):= \sum_n u_n(t)i_n(t) = \sum_b u_b(t)i_b(t),
\end{equation}
where $n$ refers to the phase number, it should be emphasized that for nonlinear and/or time-varying loads the sole interpretation of reactive power as the rate of reactive energy loses its validity. Preliminary studies show that this can be resolved by extending the reactive power balance (\ref{eq:PBreactive_st}) to include the characteristics of nonlinear and/or time-varying resistors. In case the load contains nonlinear and/or time-varying inductors and capacitors, the magnetic and electric energies need to be extended accordingly.  

\newpage

\section*{Acknowledgements}
Gerald Kaiser was supported by AFOSR Grant \#FA9550-12-1-0122 while the ideas in \cite{Kaiser} were developed.


\begin{thebibliography}{1}

\bibitem{Steinmetz1892}C.~P.~Steinmetz, ``Findet eine Phasenverschiebung im Wechselstromlichtbogen statt?,'' \em Elektrotechnische Zeitschrift \rm (1892), Issue 42, pp. 567--568. 
An English translation is available at \url{http://arxiv.org/abs/1602.06868}.

\bibitem{Pasko2012}
G.~Benysek and M.~{Pasko (editors)}.
\newblock {\em Power Theories for Improved Power Quality}.
\newblock Springer-Verlag London, 2012.

\bibitem{Budeanu1927} C.~I. Budeanu.
\newblock {\em Puissances r\'eactives et fictives}.
\newblock Inst. Romain de l'Energie, Bucharest, Romania, 1927.

\bibitem{Carozzi2005} T.~D. Carozzi, J.~E.~S Bergman and R.~L. Karlsson, ``Complex Poynting Theorem as a conservation law.'' Preprint, 2005.

\bibitem{Czarnecki1987}
L.~S. Czarnecki, ``What is wrong with the {B}udeanu concept of reactive and distortion
  power and why it should be abandoned,''
\newblock {\em IEEE Trans. Instr. Meas.}, 36(3):834--837, 1987.

\bibitem{Czarnecki1991}
L.~S. Czarnecki, ``Scattered and reactive current, voltage, and power in
  circuits with nonsinusoidal waveforms and their compensation,'' {\em IEEE
  Trans. Instr. Meas.}, vol.~40, no.~3, pp.~563--567, 1991.

\bibitem{Czarnecki2008}
L.~S. Czarnecki, ``Currents' physical components ({CPC}) concept: a fundamental
  power theory,'' {\em Przeglad Elektrotechniczny}, vol.~R84, no.~6,
  pp.~28--37, 2008.

\bibitem{Desoer}
C.~A. Desoer and E.~S. Kuh.
\newblock {\em Basic Circuit Theory}.
\newblock McGraw-Hill Book company, 1969.

\bibitem{EmanuelBook}
A.~E. Emanuel.
\newblock {\em Power Definitions and the Physical Mechanism of Power Flow}.
\newblock Wiley-IEEE Press, 2010.

\bibitem{Hanoch2006}
H.~Lev-Ari and A.~M. Stankovi\'c, ``A decomposition of apparent power in
  polyphase unbalanced networks in nonsinusoidal operation,'' {\em IEEE Trans.
  Power Sys.}, vol.~21, no.~1, pp.~438--440, 2006.

\bibitem{JeltsemaPrzeglad2016} D. Jeltsema, ``Budeanu's concept of reactive and distortion power revisited,'' {\em Przeglad Elektrotechniczny}, vol.~R92, no.~4, pp.~68--73, 2016. 

\bibitem{Kaiser} G.~Kaiser, ``Completing the complex Poynting theorem: Conservation of reactive energy in reactive time,'' 2014. \url{http://arxiv.org/abs/1412.3850}. See also
``Conservation of reactive EM energy in reactive time,'' in proc. IEEE AP-S Symposium on Antennas \& Propagation and URSI CNC/USNC Joint Meeting, Vancouver, BC, July 2015.
\url{http://arxiv.org/abs/1501.01005}.

\bibitem{Penfield}
P.~{Penfield Jr.}, R.~Spence, and S.~Duinker.
\newblock {\em Tellegen's Theorem and Electrical Networks}.
\newblock Research Monograph No. 58, MIT Press, 1970.

\bibitem{Filipski1993} P.~S.~Filipski, ``Apparent power -- a misleading quantity in the non-sinusoidal power theory: Are all non-sinusoidal power theories doomed to fail?,'' ETEP Vol. 3, No 1, Jan/Feb 1993.

\bibitem{VakmanBook}
D.~Vakman.
\newblock {\em Signals, Oscillations, and Waves. A Modern Approach}.
\newblock Artech House Inc., 1998.

\end{thebibliography}
\end{document}